\documentstyle{article}
\newcounter{eqnletter}[equation]
\catcode`\@=11
\@addtoreset{equation}{section}   
\catcode`\@=12

\begin{document}

\begin{center}

{\LARGE \bf Returns to origin of a one-dimensional random walk visiting each site an even
number of times. }
\vskip 1cm
{\large {\bf G.M. Cicuta , M.Contedini} }
\vskip 0.1 cm
Dipartimento di Fisica, Universita di Parma,\\
and INFN, Gruppo di Parma collegato alla Sezione di Milano\\
Viale delle Scienze,  43100 Parma, Italy \footnote{E-mail
addresses: cicuta@fis.unipr.it \, , \, contedini@fis.unipr.it} \\[0.5cm]
\vskip 0.1 cm
\vspace{1 cm}
\end{center}


\begin{center}
{\large {\bf Abstract }}
\end{center}
The class of random walks in one dimension, returning to the origin, restricted
by the requirement that any site visited (different from the origin) is visited an
even number of times, is analyzed in the present note. We call this class the
{\bf even-visiting random walks} and provide a closed expression to evaluate
them.

\newpage

\section{Introduction}
The number of  random walks in one dimension that originate at a given point, we may call
the origin, and after $2n$ random steps of unit lenght to the right or to the left, return
at the origin (not necessarely for the first time) is  $( 2 n )!/(n!)^2 .$  However, if we select among them
the walks where each site different from the origin is visited an even number of times,
the walks have to consist of a number of steps multiple of $4$ and their number is
more limited. Let us call such random walks {\bf the even-visiting
walks}. In these walks the origin is visited an odd number of times, if we count
both the end points of the random walk.
The purpose of this note is the description of our evaluation of
the number of the even-visiting walks
and its asymptotics ,which is given in next section.  In the rest of this section, 
we describe how these random walks are related the coefficients of
the perturbative expansion of  the resolvent of a real non-symmetric
random matrix.\\
We consider the tridiagonal matrix $M$ of order $N$ of the form
\begin{eqnarray}
M \,=\,  \left(  \begin{array}{ccccccccc}
   0 & x_1 & 0 & 0 & 0 & 0 &  \cdots & 0 & 1 \\
   1 & 0 & x_2 & 0 & 0 & 0  & \cdots & 0 & 0\\
   0 & 1 & 0 & x_3 & 0 & 0 &  \cdots & 0 & 0 \\
   0 & 0 & 1 & 0 & x_4 & 0 &  \cdots & 0 & 0 \\
   0 & 0 & 0 & 1& 0 & x_5 &  \cdots & 0 & 0\\
 \cdots & \cdots &  \cdots &  \cdots &  \cdots &  \cdots &  \cdots &  \cdots &  \cdots \\
   0 & 0 & 0 & 0 & 0 & 0 &  \cdots & 0 & x_{N-1} \\   
 x_{N} & 0 & 0 & 0 & 0 & 0 &   \cdots  & 1  & 0\\   
\end{array}   \right)  
\label{a.1}
\end{eqnarray}
The entries $x_i$ , $i=1, 2,.., N$ ,  are independent random variables with the same
probability distribution $P(x)=\frac{1}{2} \left( \delta(x-1)+\delta(x+1) \right)$. The
evaluation of the distribution of complex eigenvalues of the matrix $M$ in the
large $N$ limit is a non trivial problem, suggested to us by prof. A.Zee.
A numerical exploration of such distribution of eigenvalues , for $N=400$
 was presented in \cite{FZ}, in the context of investigations of
 one dimensional Schrodinger equation on a chain, where the  hopping
term may have a random amplitude.\\
The resolvent $G(z)$ of the random matrix $M$ is usually defined
\begin{eqnarray}
G(z)=  \frac{1}{N} {\rm Tr} <\frac{1}{z-M}>
\label{a.2}
\end{eqnarray}
where the expectation value is the mean value with the independent identically
distributed random variables $x_i$ :
\begin{eqnarray}
< f(M)>= \int \left[ \prod_{i=1}^N P(x_i) \, dx_i \right] f(M)
\label{a.3}
\end{eqnarray}
The formal perturbative expansion of $G(z)$ is
\begin{eqnarray}
G(z)= \frac{1}{N}  \sum_{k=0}^{\infty}  \frac{ <{\rm Tr} M^k >}{z^{k+1}}
\label{a.4}
\end{eqnarray}
If the order $N$ of the matrix is greater than $k$ , as we shall always
suppose , since we are interested in the limit $N \to \infty$, any term
on the diagonal of $<M^k>_{r r}$ has the same value, then the trace
merely cancels the $1/N$ factor. Let us consider the term $k=8$
\begin{eqnarray}
\sum_{a,b,c,d,e,f,g} < M_{ra} M_{ab} M_{bc} M_{cd} M_{de} M_{ef} M_{fg} M_{gr}>
\label{a.5}
\end{eqnarray}
By recalling that the non vanishing matrix elements are
 $ M_{i j}=1 $ if $j=i-1 ,$   $ M_{i j}=x_i $ if $j=i+1$ , each term in the sum (\ref {a.5})
corresponds to a walk of 8 steps, originating and ending at site $r$, with
4 steps up and 4 steps down. For instance, the sequence indicated in Fig.1,  \\
$M_{r,r+1}
M_{r+1,r+2} M_{r+2,r+1} M_{r+1,r+2} M_{r+2,r+1} M_{r+1,r} M_{r,r+1} M_{r+1,r}$ =
$x_r \cdot x_{r+1} \cdot 1 \cdot x_{r+1} \cdot 1 \cdot 1 \cdot x_r \cdot 1 = 
1$ while the
sequence in Fig.2 is $x_r \cdot x_{r+1} \cdot 1  \cdot 1 \cdot x_r  \cdot 1 \cdot x_r
 \cdot 1 $ = $x_r \cdot x_{r+1}$. Each walk corresponding to a product of
random variables $\prod_j (x_j)^{n_j}$ where all the powers $n_j$ are even
numbers, yields a contribution $+1$ , while the walks where one or several random
variables occurr with odd exponent are averaged to zero. Then $\frac{1}{N}
<{\rm Tr} M^k> $ is equal to the number of even visiting walks from a fixed site
$r$ to the same site $r$ , of $k$ steps. Because the number of steps up has to be
even, the total number of steps $k$ is multiple of $4$ and we rewrite eq.(\ref{a.4})
as
\begin{eqnarray}
G(z)= \sum_{k=0}^{\infty} \frac{c_k}{z^{4 k +1}} \quad ; \quad
c_k= \lim_{N \to \infty} \frac{1}{N} <{\rm Tr} \, M^{4 k} >
\label{a.6}
\end{eqnarray}

\section{Evaluation of the number of even visiting walks.}

Let $N(2n_r, 2 n_{r+1} , \dots , 2 n_{r+j})$ be the number of even visiting walks
corresponding to the multiplicity of the product $(x_r)^{2n_r} (x_{r+1})^{2 n_{r+1}}
\dots (x_{r+j})^{2 n_{r+j} }$. Each walk in this set visits only sites $s \geq r $
, the length of the walk is $l=4 (n_r+n_{r+1}+\dots +n_{r+j})$. The "maximum
site" visited is the site $r+j+1$ , visited $2n_{r+j}$ times ; the "minimum site"
 visited is the site $r$ ,  visited $2n_r +1$ times. The number
$N(2n_r, 2 n_{r+1} , \dots , 2 n_{r+j} , 2 n_{r+j+1} )$ is related to the
previous one $N(2n_r, 2 n_{r+1} , \dots , 2 n_{r+j})$ as it follows : new walks of
length two corresponding to $(x_{r+j+1} \cdot 1)$ may be inserted in each of the
maxima of the previous walk. Since $2 n_{r+j+1}$ identical objects are placed in
$2 n_{r+j}$ places in $ \left( 2 n_{r+j+1}+2 n_{r+j}-1 \atop 2 n_{r+j+1} \right)$
ways, we obtain
\begin{eqnarray}
&& N(2n_r, 2 n_{r+1} , \dots , 2 n_{r+j} , 2 n_{r+j+1} )=
\nonumber \\
&=&
 \left( 2 n_{r+j+1}+2 n_{r+j}-1 \atop 2 n_{r+j+1} \right)
N(2n_r, 2 n_{r+1} , \dots , 2 n_{r+j})
\label{b.1}
\end{eqnarray}
By iterating eq.(\ref{b.1}) with the initial condition $N(2 n_r)=1$ one
obtains
\begin{eqnarray}
N(2n_r, 2 n_{r+1} , \dots , 2 n_{r+j}) = \prod_{i=1}^{j-1}
 \left( 2 n_{r+i+1}+2 n_{r+i}-1  \atop 2 n_{r+i+1} \right)
\label{b.2}
\end{eqnarray}
In the same way the number $N(2 n_{r-1}, 2 n_r, 2 n_{r+1} , \dots , 2 n_{r+j})$
may be evaluated from $N(2n_r, 2 n_{r+1} , \dots , 2 n_{r+j})$ . Here the walks
 of the first number are obtained by
inserting $2 n_{r-1}$ 
walks of length two $(1 \cdot x_{r-1})$ in each of the
 $ 2 n_r +1$  minima of
the walks of the second number, obtaining
\begin{eqnarray}
&& N(2 n_{r-1}, 2 n_r, 2 n_{r+1} , \dots , 2 n_{r+j})=
\nonumber \\
&=&  \left( 2 n_{r-1}+2 n_r \atop 2 n_{r-1} \right)
N(2n_r, 2 n_{r+1} , \dots , 2 n_{r+j} )
\label{b.3}
\end{eqnarray}
Each walk contributing to $N(2 n_{r-1}, 2 n_r, 2 n_{r+1} , \dots , 2 n_{r+j})$
visits $2 n_{r-1}$ times
the "minimum site" $r-1$ , then the iteration
of the procedure leads to
\begin{eqnarray}
&& N(2 n_{r-s} , 2 n_{r-s+1}, \cdots, 2 n_{r-1}, 2 n_r, 2 n_{r+1} , \dots , 2 n_{r+j})=
\nonumber \\ &=&  \left[
\prod_{p=0}^{s-2} \left(  2 n_{r-s+p}+2 n_{r-s+p+1} -1 \atop 2 n_{r-s+p} \right)
\right] \left( 2 n_{r-1}+2 n_r \atop 2 n_{r-1} \right)
 \left[\prod_{i=1}^{j-1}
 \left( 2 n_{r+i+1}+2 n_{r+i}-1  \atop 2 n_{r+i+1} \right) \right]
\nonumber \\
\label{b.4}
\end{eqnarray}

The coefficient $c_k$ , we wish to evaluate, is the sum of several multiplicities
$N(2 n_{r-s} , 2 n_{r-s+1}, \cdots, 2 n_{r-1}, 2 n_r, 2 n_{r+1} , \dots , 2 n_{r+j})$
given above , where  $ k=2 n_{r-s} + 2 n_{r-s+1}+ \cdots +2 n_{r+j} $
corresponds to even visiting walks of $4 k $ steps. 

The evaluation may be someway simplified, by considering 
walks of fixed width, that is the difference between the "maximum site"
visited and the "minimum site" visited .
We consider the set of ordered partitions
of $k$ into positive integers $ [n_1, n_2, \cdots , n_t]$ where $k=\sum n_p$.
To each ordered sequence , or more properly to each composition ,
 $ [n_1, n_2, \cdots , n_t]$  corresponds $t+1$ classes
of even visiting walks, which are associated to the products
\begin{eqnarray}
&& (x_{r-t})^{2 n_1} (x_{r-t+1})^{2 n_2} \cdots  (x_{r-1})^{2 n_t} \quad ; 
\nonumber \\
&& (x_{r-t+1})^{2 n_1} (x_{r-t+2})^{2 n_2} \cdots  (x_{r})^{2 n_t} \quad ; 
\nonumber \\
&& \cdots \quad \cdots \quad  \cdots \quad   \cdots \quad ; \nonumber \\
&& (x_{r})^{2 n_1} (x_{r+1})^{2 n_2} \cdots  (x_{r+t-1})^{2 n_t} 
\label{b.5}
\end{eqnarray}
All walks in eq.(\ref{b.5}) have the same width $w=t$. Their multiplicities,
given in eq.(\ref{b.4}) are simply related and their sum is
\begin{eqnarray}
S_{ [n_1, n_2, \cdots , n_t ] }= \frac{ 2k}{n_1}  \prod_{i=1}^{j-1}
 \left( 2 n_{r+i+1}+2 n_{r+i}-1  \atop 2 n_{r+i+1} \right)
\label{b.6}
\end{eqnarray}

Next we sum over the ordered partitions that correspond to the 
 permutations of the
positive integers $ [n_1, n_2, \cdots , n_t]$ , finally over the different widths,
from $1$ to $k$, that is $c_k$ is the sum over the $2^{k-1}$ compositions :
\begin{eqnarray}
c_k=\sum_{t=1}^k  \sum_{ {\rm perm.}} S_{[n_1, n_2, \cdots , n_t ]} =
 \sum_{ {\rm comp.}}  S_{[n_1, n_2, \cdots , n_t ]}
\label{b.7}
\end{eqnarray}
The evaluation of eq.(\ref{b.7}) may be automated and we find
 the first coefficients $c_k$  :
\begin{eqnarray}
\begin{tabular}{||r | r ||} \hline
$c_0$ & 1       \\  \hline
$ c_1$& 2       \\  \hline
$ c_2 $ & 14      \\  \hline
$c_3$ & 116    \\  \hline
$c_4$ & 1 \,110  \\  \hline
$c_5$ & 11 \, 372 \\  \hline
$c_6$ & 123 \, 020 \\  \hline
$c_7$ & 1 \, 384 \, 168 \\  \hline
$c_8$ & 16 \, 058 \, 982 \\  \hline
$c_9$ & 190 \, 948 \, 796 \\  \hline
$c_{10}$ & 2 \,317 \, 085 \, 924 \\  \hline
$c_{11}$ & 28 \, 602 \, 719 \, 576 \\  \hline
$c_{12}$ & 358 \, 298 \, 116 \, 092 \\  \hline
$c_{13}$ & 4 \, 545 \, 807 \, 497 \, 272 \\  \hline
$c_{14}$ & 58 \, 321 \, 701 \, 832 \, 408 \\  \hline
$c_{15}$ & 755 \, 700 \, 271 \, 652 \, 816  \\  \hline
$c_{16}$ & 9 \, 878 \, 971 \, 460 \, 641 \, 414 \\  \hline
\end{tabular}
\label{b.8}
\end{eqnarray}
The ratios $c_{n}/c_{n-1}$ rise monotonically with a rate slower at higher
values of $n$. Dr.L.Molinari provided us a proof \cite{mol} that the eigenvalues 
of the random matrix $M$ , (\ref{a.1}),  are inside the square 
 with vertices $(\pm 2,0) $ , $ (0, \pm 2i) $ . This implies that $c_n \sim 16^n$. 
Let us represent the asymptotic behaviour of $r_N$ , the number of returns to the
origin of unrestricted random walks of $N$ steps,
$r_N \sim A \, \mu^N \, N^{\gamma-1} \, ,$  then $\mu=2 $ ,  $\gamma=1/2$. The
previous asymptotic behaviour for $c_n$  implies that $\mu=2 $  also for the
returns to the origin of even-visiting random walks. \\

{ \bf \large Figure Captions} \\

Fig.1 One of the even-visiting random walks, returning to site $r$ after $8$ steps, with width $w=2$.\\

Fig.2 One random walk not belonging to the class of even-visiting random walks.


\begin{thebibliography}{99}
\bibitem{mol} We plan to provide this proof together with the analysis of the
distribution of the eigenvalues in the near future.
\end{thebibliography}
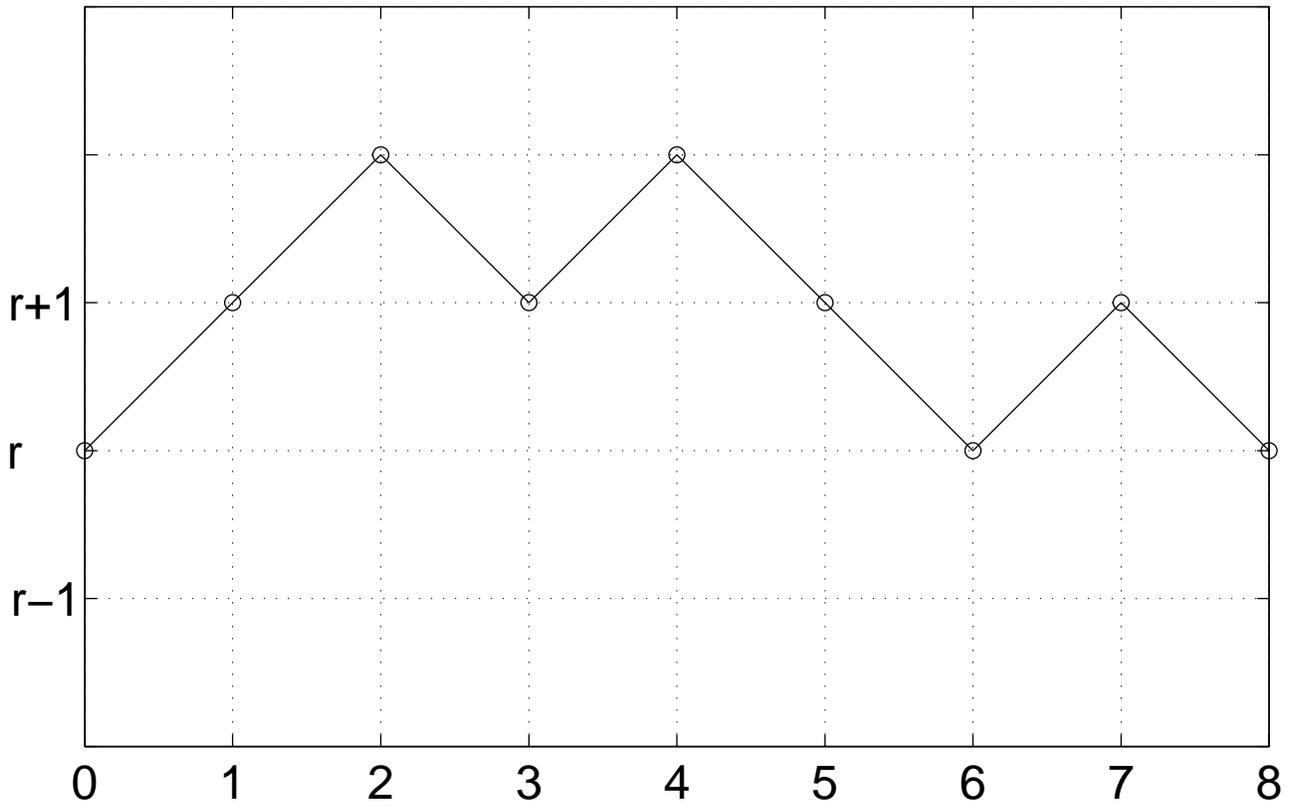
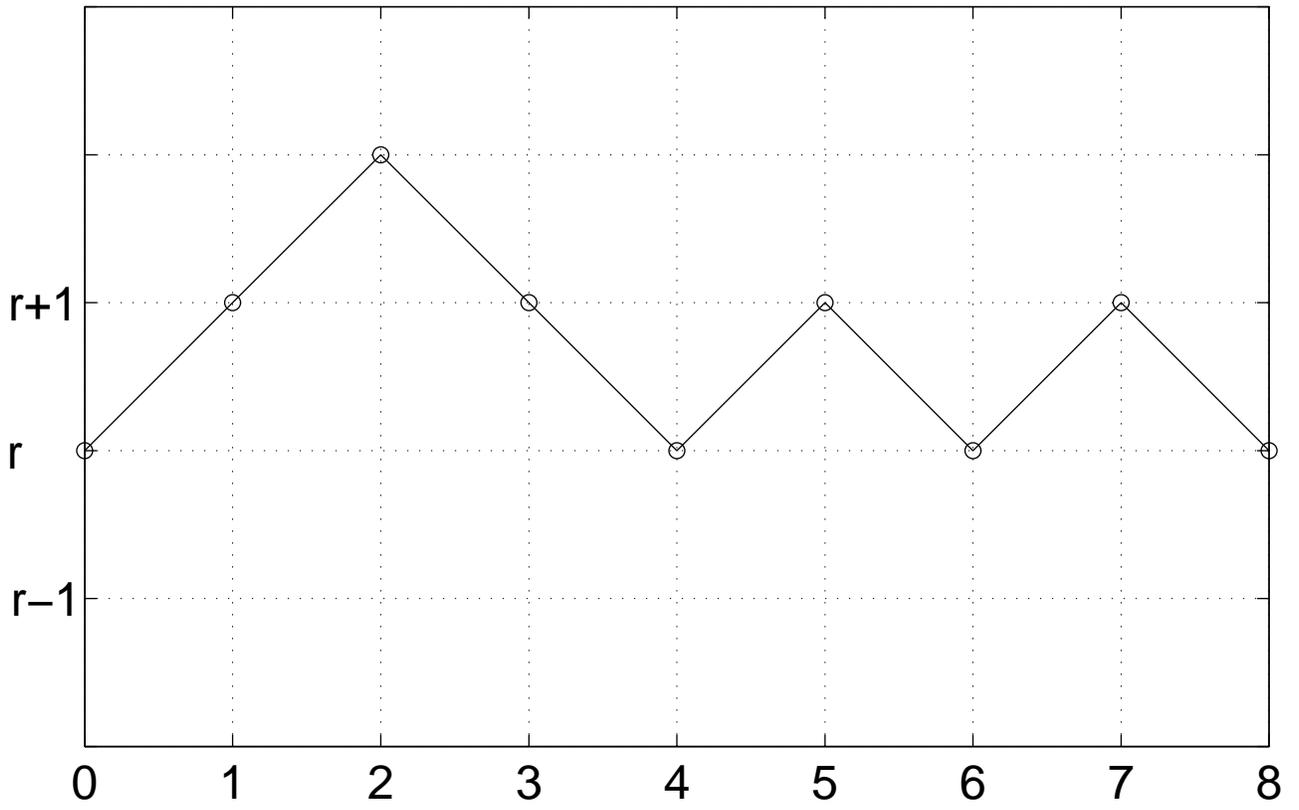
\end{document}